\documentclass[12pt]{article}
\usepackage{latexsym}
\usepackage{amsmath}
\usepackage{amssymb}
\usepackage{dsfont}
\usepackage{amsmath}
\usepackage{fancyhdr}
\usepackage{mathrsfs}
\usepackage{times}
\usepackage{graphicx}
\usepackage{theorem}

\def \C {\mathbb C}
\def \Z {\mathbb Z}
\def \R {\mathbb R}

\def \U {\mathcal U}

\def \W {\mathcal W}

\def \W {\mathcal W}

\def \p {\partial}

\begin{document}

\title{A note on K\"ahler-Ricci soliton}

\author{Xiuxiong Chen, Song Sun, Gang Tian}
\date{}
\maketitle

\begin{abstract}
 In this note we provide a proof of the following: Any
compact KRS with positive bisectional curvature is biholomorphic
to the complex projective space. As a corollary, we obtain an
alternative proof of the Frankel conjecture by using the
K\"ahler-Ricci flow.
\end{abstract}

The purpose of this note is to give a proof of the following
theorem, which does not rely on
the previous solutions of Frankel conjecture:\\

\textbf{Theorem 1}. \emph{An $n$ dimensional compact complex
manifold admitting a K\"ahler-Ricci soliton with positive
bisectional
curvature is biholomorphic to the complex projective space $\C P^n$.} \\

\textbf{Remark 2}. Since the Futaki invariant of $\C P^{n}$
vanishes, the K\"ahler-Ricci soliton must be K\"ahler-Einstein. In
addition, by a theorem of Berger (cf. \cite{B}, \cite{CT2}), it is
actually a constant multiple of the standard Fubini-Study metric.
It also follows from the uniqueness theorem in \cite{TZ1}.
\\

\textbf{Remark 3}. Using a method which originated in \cite{M},
the above theorem was proved in \cite{CT3} without using the
uniformization theorem. The proof given here is different
and we use some Morse theory.\\

As a by-product, one can use the method of K\"ahler-Ricci flow to
prove the following Frankel conjecture.\\

 \textbf{Corollary 4}. \emph{Every compact
K\"ahler manifold with positive bisectional curvature is
biholomorphic to the complex projective space.} \\

\textbf{Remark 5}. The Frankel conjecture was proved by Siu-Yau
(\cite{SY}) using harmonic maps and by Mori (\cite{M}) via
algebraic methods. In Siu-Yau's proof, by using the theorem of
Kobayashi-Ochiai (\cite{KO}), the key thing is to show the
existence of a rational curve representing the generator of
$H_2(M;\Z)/Tor$. They first proved using the second variation
formula that in the case of positive bisectional curvature, a
stable harmonic map from the sphere is either holomorphic or
anti-holomorphic. Such a harmonic map can be constructed by an
energy-minimizing process and applying Sacks-Uhlenbeck's
blowing-up analysis. Consequently, the generator of the second
homotopy class can be represented by a single holomorphic sphere.
In Mori's proof, such a rational curve was constructed using
deformation theory of curves and some
algebraic geometry of positive characteristic.\\

There has been much work on attempting to prove the Frankel
conjecture using the method of K\"ahler-Ricci flow. By a theorem
of
 Berger (c.f. \cite{B}, \cite{CT2}), it suffices to show that the flow
 converges to a K\"ahler-Einstein metric with positive bisectional
curvature. The one dimensional case was completely settled by
 R. Hamilton (\cite{H2}), B. Chow (\cite{Chow1}) and
 Chen-Lu-Tian (\cite{CLT}).
In higher dimensions, assuming $c_1>0$ and the existence of a
K\"ahler-Einstein metric, Chen-Tian (c.f. \cite{CT1}, \cite{CT2})
proved that if the initial metric has positive bisectional
curvature, then the flow converges exponentially fast to a
K\"ahler-Einstein metric with positive bisectional curvature. It
follows that the space of K\"ahler metrics with positive
bisectional curvature is path-connected. In a talk at MIT, G.
Perelman showed uniform estimates on the scalar curvature and
diameter along the K\"ahler-Ricci flow and announced that the flow
with arbitrary initial metric converges to a K\"ahler-Einstein
metric if there exists one. In \cite{TZ2}, using an estimate of
Perelman, Tian-Zhu gave a proof of the convergence of the
K\"ahler-Ricci flow on any K\"ahler-Einstein manifold via a
different method. In \cite{C1}, X. Chen proved that an irreducible
Fano K\"ahler manifold with positive orthogonal bisectional
curvature is biholomorphic to $\C P^n$ using the Frankel
conjecture. Recently, Phong-Song-Sturm-Weinkove (\cite{PSSW})
proved that the K\"ahler-Ricci flow starting from a metric of
positive bisectional
curvature will converge to a K\"ahler-Einstein metric under various extra conditions. \\

Our proof of theorem 1 uses induction on the dimension $n$,
starting from the trivial case $n=0$. Denote by $A_n$ the
statement of theorem 1 for dimension $n$, and by $B_n$ the
statement of corollary 4 for dimension $n$. We will prove first
that $A_n$ implies $B_n$, and then that $B_{k}(k<n)$ implies
$A_n$. So theorem 1 and corollary 4 are proved simultaneously.\\

\emph{Proof of Corollary 4($A_n\rightarrow B_n$)}. Suppose that
$(M, J, g)$ is a $n$ dimensional K\"ahler manifold with positive
bisectional curvature. Then the tangent bundle $TM$ is ample, so
in particular, $c_1(TM)>0$. By applying a Bochner type formula
(c.f. \cite{BG}), we know $b_2(M)=1$. In fact, $b_{2k}=1$ and
$b_{2k+1}=0$ for all $k=0,\cdots, n$. Therefore we can assume the
K\"ahler form $\omega$ lies in the canonical class $c_1(TM)$. Now
we run the K\"ahler-Ricci flow:
$$\frac{\partial}{\partial t}g_{i\bar{j}}(t)=g_{i\bar{j}}(t)-R_{i\bar{j}}(t),$$
with $g(0)=g$. In \cite{Cao}, H. D. Cao proved that the flow
exists globally. By S. Bando (\cite{Ba}) in dimension 3 and N. Mok
(\cite{Mok}) in all dimensions, the positivity of the bisectional
curvature is preserved under the K\"ahler-Ricci flow. To study
convergence, we make use of a theorem of Perelman that the scalar
curvature and the diameter are uniformly bounded along the flow
(cf. \cite{ST}). Thus the bisectional curvature is also uniformly
bounded along the flow, then all higher order derivative estimates
of the curvature follow easily (cf. \cite{CT1}). Hence there exist
a subsequence $g(t_i)$ and diffeomorphisms $f_i$, such that
$f_i^*g(t_i)$ converges to $g_{\infty}$ smoothly. Moreover, we can
assume $f_i^*J$ converges smoothly to $J_{\infty}$(possibly
different from $J$) and $g_{\infty}$ is K\"ahler with respect to
$J_{\infty}$. Using Perelman's $\W$ functional, we can prove that
$(g_{\infty}, J_{\infty})$ is a K\"ahler-Ricci soliton (cf.
\cite{Se}). Clearly it has non-negative bisectional curvature. We
claim it actually has positive bisectional curvature. Indeed, by
Lemma 6 proved near the end of this note, it suffices to show that
the Ricci curvature of $g_{\infty}$ is positive. It follows from
the strong maximum principle for tensors along the Ricci flow (cf.
\cite{H1}) that if $Ric(g_{\infty})$ has a null direction at some
point, then the manifold $(M, g_{\infty}, J_{\infty})$ splits
holomorphically isometrically into a product of a flat factor $P$
and another factor $Q$ with strictly positive Ricci curvature.
Both $P$ and $Q$ are K\"ahler manifolds, so $b_2(P), b_2(Q)\geq
1$. Hence $b_2(M)\geq 2$, which is a contradiction. Therefore we
have proved the claim. Now by Theorem 1, we know $(M, J_{\infty})$
is biholomorphic to $\C P^n$. Since $\C P^n$ has trivial local
deformation, i.e. $H^1(\C P^n,\Theta )=0$ (c.f. \cite{Bott},
\cite{Kod}), $(M, J)$ is also biholomorphic to $\C P^n$, and this
proves Corollary 4. Note that by \cite{CT2}, the K\"ahler Ricci
flow here indeed converges exponentially fast to the Fubini-Study
metric.
$\square$\\\\

Now we come to prove that $B_k(k<n)$ implies $A_n$. It uses the
theorem of Kobayahsi-Ochiai as in the previous solutions for the
Frankel conjecture. By Bishop-Goldberg (c.f. \cite{BG}), a compact
K\"ahler manifold with positive bisectional curvature has the
second Betti number equal to $1$. Suppose an $n$ dimensional
compact complex manifold $(M, J)$ admits a K\"ahler-Ricci soliton
$g$, with the K\"ahler form $\omega(\cdot,
\cdot)=g(J\cdot,\cdot)$, i.e.
$$Ric(\omega)=\omega+\sqrt{-1}\partial\overline{\partial}f,$$
where $f$ is a real function whose gradient is a holomorphic
vector field, i.e. the (2,0) part of $\nabla\nabla f$ vanishes. In
particular, we know $f$ is a  Morse-Bott function (see \cite{F}),
i.e. the set of critical points of $f$ consists of smooth compact
submanifolds of $M$ on which $Hess(f)$ is non-degenerate along the
normal directions. We may assume $f$ is not a constant, since
otherwise $(M, J, g)$ is K\"ahler-Einstein and $(M, J)$ is
biholomorphic to the complex projective space by \cite{B}. A
critical point of $f$ is the same as a zero of the vector field
$Im \nabla f$. Since $Im \nabla f$ is Killing, the critical
submanifolds are totally geodesic in $M$; since $Im \nabla f$ is
also holomorphic, the critical submanifolds
 are all K\"ahler
submanifolds, and the Morse indices(the number of negative
eigenvalues of $Hess(f)$) are all even. By Kobayashi-Ochiai
\cite{KO} (c.f. \cite{SY}), to prove theorem 1, it suffices to
show the generator of $\pi_2(M)\simeq H_2(M;\Z)$ modulo torsion
could be represented by a rational curve, i.e. a holomorphic map
from $\C P^1$ to $M$. Obviously the critical submanifolds are all
K\"ahler manifolds with positive bisectional curvature with
dimension less than $n$, thus by induction hypothesis they are all
biholomorphic to the complex projective spaces of various
dimension. By a general theorem of Frankel(see \cite{F}), the
following holds:
$$b_i(M)=\sum_{\alpha}b_{i-\lambda_{\alpha}}(F_{\alpha}),$$
where $F_{\alpha}$'s are all critical submanifolds of $f$, and
$\lambda_{\alpha}$ is the Morse index of $F_{\alpha}$. Let $i=0,
2$, we have
$$1=b_0(M)=\sum_{\alpha:\lambda_{\alpha}=0}b_0(F_{\alpha}),$$ and
$$1=b_2(M)=\sum_{\alpha:\lambda_{\alpha}=0}b_{2}(F_{\alpha})+\sum_{\alpha:\lambda_{\alpha}=2}b_{0}(F_{\alpha}).$$

 Therefore we
have two alternatives: either there are no critical submanifolds
of index $2$, or there is exactly one critical submanifold of
index $2$. In both cases the minimal submanifold of $f$ (i.e. the
set of critical points achieving the minimum of $f$) is connected.
In the first case, the minimal submanifold is a $\C P^k$ for some
$k\geq 1$, and the complement of the stable manifold\footnote{The
stable manifold of a critical submanifold $N$ is defined to be the
set of points which will converge to $N$ under the negative
gradient flow of $f$. It is an open manifold of dimension
complementary to the Morse index of $N$.} of the minimal
submanifold has codimension at least 4. Hence a generic
representative of the generator of $\pi_2(M)$(modulo torsion) will
be deformed to a surface on the minimal manifold under the
negative gradient flow of $f$.  It is easy to see the latter must
represent the generator of $\pi_2(\C P^k)$ and is thus homotopic
to a rational curve. In the second case, the minimal submanifold
must be a point. Denote the unique critical submanifold of index
$2$ by $M_2$, and the minimal submanifold by $M_0$. Then all the
critical submanifolds other than $M_2$ and $M_0$ have Morse
indices at least 4.
Denote
$Stab(N)$ the stable manifold of the critical submanifold $N$.
Then we have the following stratification:
$$M=M_0 \cup M_2\cup_{\alpha}Stab(N_{\alpha}),$$ where $N_{\alpha}$ runs
over all critical submanifolds of $M$ with Morse index greater or
equal to $4$. We claim that the closure $\overline{Stab(M_2)}$ is
an analytic subvariety of $M$. Indeed, $Stab(M_2)$ is a complex
submanifold of $M\setminus\cup_{\alpha}Stab(N_{\alpha})$ with
complex dimension $n-1$(Using the holomorphic gradient flow of $f$
it suffices to prove this near $M_2$. The latter follows from the
fact that $X=\nabla^{1,0} f$ is linearizable at its critical
points (c.f. \cite{Bryant}). More precisely, near a point $p$ in
$M_2$, by \cite{Bryant}, we can find a coordinate chart
$(\U,\{z^i\})$ such that $X=a_i z^i\frac{\partial}{\partial z^i}$,
where $a_i\geq 0$, $i=1, \cdots, n-1$ and $a_n<0$.Thus
$Stab(M_2)\cap\U$ is given by the hyperplane $z_n=0$). Since
$\cup_{\alpha}Stab(N_{\alpha})$ is the union of finitely many
complex subvarieties of complex codimension at least $2$ in $M$,
the analyticity of $\overline{Stab(M_2)}$ follows from the Levi
extension theorem(see \cite{GH}). Therefore, the cycle
$\overline{Stab(M_2)}$ represents an element in $H_{2n-2}(M; \Z)$,
and it is smooth near $M_2$.
 On the other hand,%
near a point $p$ in $M_2$, we have the afore-mentioned linearized
coordinate $(\U,\{z^i\})$,  then locally the unstable manifold of
$p$ is given by the $z_n$ axis. Clearly the orbits of $Im\nabla f$
on the $z_n$ axis are all periodic. Pick such an orbit $\theta:
S^1\rightarrow M$, then we can construct a map $R:
S^1(\simeq\R/T\Z)\times \R\rightarrow M$ by defining $R(s, t)$ to
be $\phi_t(\theta(s))$, where $\phi_t$ is the integral curve of
the negative gradient flow of $f$, and $T$ is the period of the
orbit of $\theta$. Since $[\nabla f, J\nabla f]=0$, $R$ is
holomorphic. By the Riemann removable singularity theorem, $R$
extends to a rational curve, still denoted by $R: \C
P^1\rightarrow M$. It coincides with the unstable manifold of $p$
in $\U$, thus $R$ is smoothly embedded near $p$, intersecting
 $\overline{Stab(M_2)}$ transversally at exactly one point $p$. So
the intersection number of integral homology classes represented
by $\overline{Stab(M_2)}$ and $R$ is $1$. It follows that $R$
represents a generator of $\pi_2(M)$(modulo torsion). Now the
remaining follows from standard arguments. By a well-known result
of Grothendieck, the quotient bundle $R^* (TM/T\C P^1)$ splits
into a direct sum of line bundles $Q_2,\cdots, Q_n$, each $Q_i$
$(2\leq i\leq n)$ is positive. It follows that
$$c_1(R^*TM)=c_1(\C P^1)+ \sum_{i=2}^n c_1(Q_i).$$ Hence $c_1(TM)$
evaluated at $[R]$ is bigger or equal to $n + 1$. That is,
$c_1(TM)\geq(n+1)c_1(F)$, where $F$ is the line bundle on $M$ such
that $\langle c_1(F),[R]\rangle=1$. By the result of
Kobayashi-Ochiai (\cite{KO}),
$M$ is biholomorphic to $\C P^n$. This finishes the proof of Theorem 1. $\square$\\

We need the following lemma in the proof of Corollary 4.\\

\textbf{Lemma 6}. \emph{Along the K\"ahler-Ricci flow, suppose
there are
 positive constants $T_0$ and $C$, such that for $t\in [T_0, +\infty)$, $g(t)$ has positive bisectional curvature
and the Ricci curvature of $g(t)$ satisfies $Ric(g(t))\geq C\cdot
g(t)$. Then the bisectional curvature of $g(t)$ is uniformly
bounded
below from 0.}\\

\emph{Proof}. We follow an idea in \cite{C1} which is essentially
an application of the maximum principle. Along the K\"ahler-Ricci
flow, we have the following evolution equations:
\begin{equation*}
\frac{\p g_{i\bar{j}}}{\p t}=g_{i\bar{j}}-R_{i\bar{j}};
\end{equation*}
\begin{equation*}
\frac{\p R}{\p t}=\Delta R+|Ric|^2-R;
\end{equation*}
\begin{equation*}
\frac{\p Ric}{\p t}=\Delta Ric+Ric\cdot Rm-Ric^2;
\end{equation*}
\begin{eqnarray}
\label{eq:evo-curv} \frac{\p R_{i\bar{j}k\bar{l}}}{\p t}&=&\Delta
R_{i\bar{j}k\bar{l}}+R_{i\bar{j}p\bar{q}}R_{q\bar{p}k\bar{l}}-R_{i\bar{p}k\bar{q}}R_{p\bar{j}q\bar{l}}+
R_{i\bar{l}p\bar{q}}R_{q\bar{p}k\bar{j}}+R_{i\bar{j}k\bar{l}}\nonumber\\&&-\frac{1}{2}(R_{i\bar{p}}R_{p\bar{j}k\bar{l}}
+R_{p\bar{j}}R_{i\bar{p}k\bar{l}}+R_{k\bar{p}}R_{i\bar{j}p\bar{l}}+R_{p\bar{l}}R_{i\bar{j}k\bar{p}}).
\end{eqnarray}
Here we define $$(Ric\cdot
Rm)_{i\bar{j}}=R_{l\bar{k}}R_{i\bar{j}k\bar{l}},$$ and
$$(Ric^2)_{i\bar{j}}=R_{i\bar{k}}R_{k\bar{j}}.$$
Now we put $S=Rm-\mu(g*Ric)$, where $\mu$ is a function of $t$,
and
$$(g*Ric)_{i\bar{j}k\bar{l}}=g_{i\bar{j}}R_{k\bar{l}}+g_{k\bar{l}}R_{i\bar{j}}+g_{i\bar{l}}R_{k\bar{j}}
+g_{k\bar{j}}R_{i\bar{l}}.$$ Then
$$S_{k\bar{l}}=(1-(n+2)\mu)R_{k\bar{l}}-\mu R\cdot g_{k\bar{l}},$$
i.e
$$Sic=(1-(n+2)\mu)Ric-\mu R\cdot g,$$
where $Sic_{i\bar{j}}=g^{k\bar{l}}S_{i\bar{j}k\bar{l}}$.
Therefore, by a straightforward calculation, we obtain\\
\begin{eqnarray*} \frac{\p R_{i\bar{j}k\bar{l}}}{\p t}&=&
\square S_{i\bar{j}k\bar{l}} +\mu g*(\Delta
Ric)_{i\bar{j}k\bar{l}}+\mu
(g*Ric)_{i\bar{j}k\bar{l}}+\mu[(g*(Ric\cdot
 S))_{i\bar{j}k\bar{l}}+(Ric*Sic)_{i\bar{j}k\bar{l}}]+I\\\\
 &&-\mu[(Ric*Ric)_{i\bar{j}k\bar{l}}+(Ric^2*g)_{i\bar{j}k\bar{l}}],\\
\end{eqnarray*}
where $\square S_{i\bar{j}k\bar{l}}$ denotes the right side of
(\ref{eq:evo-curv}) with $R_{i\bar{j}k\bar{l}}$ replaced by
$S_{i\bar{j}k\bar{l}}$, and\\
\begin{eqnarray*}
I&=&\mu^2[(g*Ric)_{i\bar{j}p\bar{q}}(g*Ric)_{q\bar{p}k\bar{l}}-
(g*Ric)_{i\bar{p}k\bar{q}}(g*Ric)_{p\bar{j}q\bar{l}}\\\\&&+(g*Ric)_{i\bar{l}p\bar{q}}(g*Ric)_{q\bar{p}k\bar{j}}].\\
\end{eqnarray*}
Finally we calculate:\\
\begin{eqnarray*}
-\frac{\p }{\p t}(\mu(g*Ric)_{i\bar{j}k\bar{l}})&=&
-\mu'(g*Ric)_{i\bar{j}k\bar{l}}-\mu(\frac{\p g}{\p
t}*Ric)_{i\bar{j}k\bar{l}}-\mu(g*\frac{\p Ric}{\p
t})_{i\bar{j}k\bar{l}}\\\\&=&
-\mu'(g*Ric)_{i\bar{j}k\bar{l}}-\mu(g*Ric)_{i\bar{j}k\bar{l}}+\mu(Ric*Ric)_{i\bar{j}k\bar{l}}\\\\&&
-\mu(g*(\Delta Ric))_{i\bar{j}k\bar{l}}-\mu(g*(Ric\cdot
Rm))_{i\bar{j}k\bar{l}}+\mu(g*(Ric^2))_{i\bar{j}k\bar{l}}.
\end{eqnarray*}\\
It follows that\\
\begin{eqnarray*}
\frac{\p S_{i\bar{j}k\bar{l}}}{\p t}&=& \square
S_{i\bar{j}k\bar{l}}+\mu g*(\Delta Ric)_{i\bar{j}k\bar{l}}+\mu
(g*Ric)_{i\bar{j}k\bar{l}}+\mu[(g*(Ric\cdot
S))_{i\bar{j}k\bar{l}}+(Ric*Sic)_{i\bar{j}k\bar{l}}]+I\\\\&&
-\mu[(Ric*Ric)_{i\bar{j}k\bar{l}}+(Ric^2*g)_{i\bar{j}k\bar{l}}]
-\mu'(g*Ric)_{i\bar{j}k\bar{l}}-\mu(g*Ric)_{i\bar{j}k\bar{l}}\\\\&&+\mu(Ric*Ric)_{i\bar{j}k\bar{l}}
-\mu(g*(\Delta Ric))_{i\bar{j}k\bar{l}}-\mu(g*(Ric\cdot
Rm))_{i\bar{j}k\bar{l}}+\mu(g*(Ric^2))_{i\bar{j}k\bar{l}}\\\\&=&
\square S_{i\bar{j}k\bar{l}}+\mu^2(g*(Ric\cdot
(g*Ric))_{i\bar{j}k\bar{l}}+\mu(Ric*Sic)_{i\bar{j}k\bar{l}}-\mu'(g*Ric)_{i\bar{j}k\bar{l}}+I\\\\&=&
\square
S_{i\bar{j}k\bar{l}}+\mu(Ric*Ric)_{i\bar{j}k\bar{l}}+\mu^2(g*(Ric\cdot
(g*Ric)))_{i\bar{j}k\bar{l}}\\\\&&-\mu^2(Ric*((n+2)Ric+R\cdot
g)_{i\bar{j}k\bar{l}}-\mu'(g*Ric)_{i\bar{j}k\bar{l}}+I.\\
\end{eqnarray*}
Notice that $I=O(\mu^2)$. Since $Ric(g(t))\geq C\cdot g(t)$ for
$t\geq T_0$, if $\mu(t)\equiv\mu(0)\equiv\mu>0$ is sufficiently
small, we have
$$(\frac{\p}{\p t}-\square)S_{i\bar{j}k\bar{l}}\geq0.$$
Now by the maximum principle as in \cite{Mok}, we see that
$S_{i\bar{i}j\bar{j}}\geq0$ for all $t\geq T_0$. Thus
$R_{i\bar{i}j\bar{j}}\geq \mu (g*Ric)_{i\bar{i}j\bar{j}}\geq
\mu C(g*g)_{i\bar{i}j\bar{j}}$. $\square$\\

 \textbf{Remark 7}. By a more careful study of critical
submanifolds of $f$ and using the localization formula, it may be
possible to give a direct proof of vanishing of the Futaki
invariant associated to the holomorphic vector field induced by
$f$. If so, one would obtain a proof of Theorem 1
which does not rely on the theorem of Kobayashi-Ochiai.\\

\noindent{\bf Acknowledgments:} Both the first and third authors
are partially supported by NSF grants. We would like to thank the
anonymous referee for pointing out an error in the earlier version
of this paper.

Xiuxiong Chen, Van Vleck Hall, University of Wisconsin-Madison,
480 Lincoln Drive, Madison, WI 53706.\\ Email:
xxchen@math.wisc.edu\\\\
 Song Sun, Van Vleck Hall, University of
Wisconsin-Madison, 480 Lincoln Dr, Madison, WI 53706.\\ Email:
ssun@math.wisc.edu\\\\
 Gang Tian, Fine Hall, Princeton University,
Princeton, NJ 08544.\\ Email: tian@math.princeton.edu\\\\
\end{document}